\documentclass[10pt]{amsart}
\usepackage{mathrsfs}
\usepackage{amsfonts} 
\usepackage{graphicx,latexsym,bm,amsmath,amssymb,verbatim,multicol,lscape}
\theoremstyle{definition}

\theoremstyle{remark}

\numberwithin{equation}{section}

\begin{document}
\title[Constructing permutation polynomials over finite fields]
{Constructing permutation polynomials over finite fields}%
\author{XIaoer Qin}
\address{Mathematical College, Sichuan University, Chengdu 610064, P.R. China
and College of Mathematics and Computer Science, Yangtze Normal University,
Chongqing 408100, P.R. China}
\email{qincn328@sina.com}
\author{Shaofang Hong$^*$}
\address{Yangtze Center of Mathematics, Sichuan University, Chengdu 610064, P.R. China
and Mathematical College, Sichuan University, Chengdu 610064, P.R. China}
\email{sfhong@scu.edu.cn, s-f.hong@tom.com, hongsf02@yahoo.com}
\thanks{$^*$Hong is the corresponding author and was supported partially by the
Ph.D. Programs Foundation of Ministry of Education of China Grant \#20100181110073}
\date{\today}%
\begin{abstract}
In this paper, we construct several new permutation polynomials over
finite fields. First, using the linearized polynomials, we construct
the permutation polynomial of the form
$\sum_{i=1}^k(L_{i}(x)+\gamma_i)h_i(B(x))$ over ${\bf F}_{q^{m}}$,
where $L_i(x)$ and $B(x)$ are linearized polynomials. This extends a
theorem of Coulter, Henderson and Matthews. Consequently, we
generalize a result of Marcos by constructing permutation
polynomials of the forms $x h(\lambda_{j}(x))$ and $x
h(\mu_{j}(x))$, where $\lambda_{j}(x)$ is the $j$-th elementary
symmetric polynomial of $x, x^{q}, ..., x^{q^{m-1}}$ and
$\mu_{j}(x)=\textup{Tr}_{{\bf F}_{q^{m}}/{\bf F}_{q}}(x^{j})$.
This answers an open problem raised by Zieve in 2010.
Finally, by using the linear translator, we construct the permutation
polynomial of the form $L_1(x)+L_{2}(\gamma)h(f(x))$ over ${\bf
F}_{q^{m}}$, which extends a result of Kyureghyan.
\end{abstract}

\maketitle

\section{\bf Introduction}
Let ${\bf F}_{q}$ denote the finite field of characteristic $p$ with
$q$ elements ($q=p^{n},n\in \textbf{N}$), and let ${\bf
F}_{q}^{*}:={\bf F}_{q}\setminus \{0\}$. Let ${\bf F}_{q}[x]$ be the
ring of polynomials over ${\bf F}_{q}$ in the indeterminate $x$. If
the polynomial $f(x)\in {\bf F}_{q}[x]$ induces a bijective map from
${\bf F}_{q}$ to itself, then $f(x)$ is called a {\it permutation
polynomial} of ${\bf F}_{q}$. Permutation polynomials have been an
interesting subject of study in the area of finite fields for many
years. Particularly, permutation polynomials have many important
applications in coding theory \cite{[L]}, cryptography \cite{[SH]}
and combinatorial design theory. Information about properties,
constructions and applications of permutation polynomials may be
found in the book of Lidl and Niederreiter \cite{[LN]}.

Let $m>1$ be a given integer. By $\textup{Tr}_{{\bf F}_{q^m}/{\bf F}_{q}}(x)$
we denote the {\it trace} from ${\bf F}_{q^m}$ to ${\bf F}_{q}$, that is
$$\textup{Tr}_{{\bf F}_{q^m}/{\bf F}_{q}}(x)=x+x^{q}+\cdots+x^{q^{m-1}}.$$
A polynomial of the form $$L(x)=\sum^{m-1}_{i=0}a_{i}x^{q^{i}}\in
{\bf F}_{q^m}[x]$$ is called a \emph{linearized polynomial} over ${\bf F}_{q^m}$.
It is well known that a linearized polynomial $L(x)$ is a permutation
polynomial of ${\bf F}_{q^m}$ if and only if the set of roots in ${\bf F}_{q^m}$
of $L(x)$ equals $\{0\}$ (see, for example, Theorem 7.9 of \cite{[LN]}).
Throughout, $L(x)$ denotes a linearized polynomial.

To find new classes of permutation polynomials is one of the open
problems raised by Lidl and Mullen in \cite{[LM]}. There has been
significant progress in finding new permutation polynomials. Wan and
Lidl \cite{[WL]}, Masuda and Zieve \cite{[MZ]} and
Zieve \cite{[Z2]} constructed permutation polynomials
of the form $x^{r}f(x^{(q-1)/d})$ and studied their group structure. Zieve
\cite{[Z]} characterized the permutation polynomial of the form
$x^r(1+x^v+x^{2v}+...+x^{kv})^t$. Ayad, Belghaba and Kihel
\cite{[ABK]} obtained some permutation binomials and proved the
bound of $p$, if $ax^n+x^m$ permutes ${\bf F}_{p}$. A number of
classes of permutation polynomials related to the trace functions
were constructed. Recently, Coulter, Henderson and Matthews
\cite{[CHM]} constructed the permutation polynomials of the form
$L(x)+xh(\textup{Tr}_{{\bf F}_{q^m}/{\bf F}_{q}}(x))$. Marcos
\cite{[M]} obtained permutation polynomials of the form
$bL(x)+\gamma h(\textup{Tr}_{{\bf F}_{q^m}/{\bf F}_{q}}(x))$. Zieve
\cite{[Z1]} presented rather more general versions of the first four
constructions from \cite{[M]}. But how to extend the $5$-th
construction from \cite{[M]} to a more general version
is an interesting open problem raised in
\cite{[Z1]}. For some other permutation polynomials constructed by
using the trace function, the readers are referred to \cite{[CK]}.

The main goal of the present paper is to construct new classes of
permutation polynomials over finite fields. In Section 2, we
construct some permutation polynomials using linearized polynomials.
In fact, we obtain a characterization so that
$\sum_{i=1}^k(L_{i}(x)+\gamma_i)h_i(B(x))\in {\bf F}_{q^m}[x]$ with
$L_i(x)$ and $B(x)$ being linearized polynomials, is a permutation
polynomial. See Theorem 2.1 below, which extends the results
obtained by Coulter, Henderson and Matthews \cite{[CHM]} and by
Marcos \cite{[M]}, respectively.

For any integer $j$ with $1\le j\le m-1$, let $\mu_{j}(x)
=\textup{Tr}_{{\bf F}_{q^m}/{\bf F}_{q}}(x^{j})$ and
$\lambda_{j}(x)=\sigma_{j}(x,x^{q},\ldots,x^{q^{m-1}})$, where
$\sigma_{j}(x,x^{q},\ldots,x^{q^{m-1}})$ is the $j$-th elementary
symmetric polynomial of $x, x^{q}, ..., x^{q^{m-1}}$. Marcos
\cite{[M]} used the function $\lambda(x)(=\lambda_{2}(x)$ or
$\mu_{2}(x))$ to construct permutation polynomials and only got some
sufficient conditions so that $xh(\lambda (x))$ to be a permutation
polynomial. In Section 3, using $\lambda_{j}(x)$ and $\mu_{j}(x)$,
we extend this result of Marcos \cite{[M]} by giving the sufficient
and necessary conditions so that $xh(\lambda_{j}(x))$ and
$xh(\mu_{j}(x))$ to be permutation polynomials. This answers an
open problem raised by Zieve in \cite{[Z1]}.

In Section 4, by using the technique of linear translator (see
Section 4 for its definition), we construct the permutation polynomial
of the form $L_1(x)+L_{2}(\gamma)h(f(x))$. This result generalizes
one of the main results in \cite{[K]}.

\section{\bf Permutation polynomials constructed by the linearized polynomials}

In this section, we construct a new class of permutation polynomials
involving linearized polynomials. We need the following known facts in the sequel.\\

\noindent{\bf Lemma 2.1.} {\it Let $B(x)\in {\bf F}_{q}[x]$ and
$L(x)\in {\bf F}_{q}[x]$ be linearized polynomials. Then for any
$a\in {\bf F}_{q}$ and
 $x$ and $y\in {\bf F}_{q^m}$, $aB(x)=B(ax)$, $B(x+y)=B(x)+B(y)$ and
 $B(L(x))=L(B(x))$.}\\

We can now give the first main result of this paper.\\

\noindent{\bf Theorem 2.1.} {\it For $ 1\leq i\leq k$, let
$\gamma_i\in {\bf F}_{q^m}$ and let $L_{i}(x)$, $B(x)\in {\bf
F}_{q}[x]$ be linearized polynomials. Let $h_i(x)\in {\bf
F}_{q^m}[x]$ be such that $h_i(B({\bf F}_{q^m}))\subseteq {\bf
F}_{q}$. Then $F(x):=\sum_{i=1}^k(L_{i}(x)+\gamma_i)h_i(B(x))$
is a permutation polynomial over ${\bf F}_{q^m}$
if and only if each of the following is true: \\
{\rm (1).} $\sum_{i=1}^k(L_{i}(x)+B(\gamma_i))h_i(x)$ permutes $B({\bf F}_{q^m})$.\\
{\rm (2).} For any $y\in B({\bf F}_{q^m})$, $\sum_{i=1}^k
L_{i}(x)h_i(y)=0$ and $B(x)=0$ with $x\in {\bf F}_{q^m}$ are both
true if and only if $x=0$. }

\begin{proof}
First we show the sufficiency part. Assume that (1) and (2) hold.
Suppose that there exist two elements $\alpha$ and $\beta \in {\bf F}_{q^m}$
such that $F(\alpha)=F(\beta)$. Thus $B(F(\alpha))=B(F(\beta))$. That is,
$$B\Big(\sum_{i=1}^k(L_{i}(\alpha)+\gamma_i)h_i(B(\alpha))\Big)=
B\Big(\sum_{i=1}^k(L_{i}(\beta)+\gamma_i)h_i(B(\beta))\Big).\eqno(2.1)$$
Then Lemma 2.1 applied to both sides of (2.1) gives us that
$$\sum_{i=1}^k(L_{i}(B(\alpha))+B(\gamma_i))h_i(B(\alpha))
=\sum_{i=1}^k(L_{i}(B(\beta))+B(\gamma_i))h_i(B(\beta)).\eqno(2.2)$$
Since $\sum_{i=1}^k(L_{i}(x)+B(\gamma_i))h_i(x)$ permutes $B({\bf
F}_{q^m})$, it follows from (2.2) that $B(\alpha)=B(\beta)$. Write
$t:=B(\alpha)=B(\beta)$. Then $t\in B({\bf F}_{q^m})$ and
$B(\alpha-\beta)=0$. Since $F(\alpha)=F(\beta)$, one has
$$\sum_{i=1}^kL_{i}(\alpha-\beta)h_i(t)=0.$$
Now applying condition (2) to $\alpha-\beta$, we obtain that
$\alpha-\beta=0$ which implies that $\alpha=\beta$. Hence $F(x)$ is
a permutation polynomial over ${\bf F}_{q^m}$. The sufficiency part is proved.

Let us now show the necessity part. Let $F(x)$ be a permutation
polynomial of ${\bf F}_{q^m}$. First we prove that (1) is true. To do so, we let
$B(x)$ act on $F(x)$ for $x\in {\bf F}_{q^m}$, and then by Lemma 2.1 we get that
$$B(F(x))=\sum_{i=1}^k(L_{i}(B(x))+B(\gamma_i))h_i(B(x)).\eqno(2.3)$$

Since $F(x)$ is a permutation polynomial of ${\bf F}_{q^m}$, we have
$$|\{B(F(x)):x\in {\bf F}_{q^m}\}|=|\{B(x):x\in {\bf F}_{q^m}\}|=|B({\bf F}_{q^m})|.\eqno(2.4)$$
Hence by (2.3) and (2.4),
$$|\{\sum_{i=1}^k(L_{i}(B(x))+B(\gamma_i))h_i(B(x)):x\in {\bf F}_{q^m}\}|
=|B({\bf F}_{q^m})|.$$ This concludes that
$\sum_{i=1}^k(L_{i}(x)+B(\gamma_i))h_i(x)$ permutes $B({\bf
F}_{q^m})$. Thus (1) is proved.

It remains to show that (2) is true. For this purpose, we
assume that for certain $y\in B({\bf F}_{q^m})$ and $x\in {\bf F}_{q^m}$,
we have $\sum_{i=1}^kL_{i}(x)h_i(y)=0$ and $B(x)=0$. We can take two elements
$\alpha$ and $\beta\in {\bf F}_{q^m}$ satisfying that $B(\alpha)=B(\beta)=y$.
Then $B(\alpha-\beta)=0$. But $B(x)=0$. Therefore $\alpha-\beta$ and
$x$ are both in the kernel ${\rm ker}(B)$ of $B(x)$. So we can write
$x=\alpha-\beta+z$ for some $z\in {\rm ker}(B)$. Since $\sum_{i=1}^k
L_{i}(x)h_i(y)=0$, we infer that
$$\sum_{i=1}^kL_{i}(\alpha-\beta+z)h_i(y)=0.\eqno(2.5)$$

On the other hand, since $z\in {\rm ker}(B)$, one has $B(z)=0$,
which implies that $B(\beta-z)=B(\alpha)=y$. It then follows
immediately that
\begin{align*}
&F(\alpha)-F(\beta-z)\\
&=\sum_{i=1}^k (L_{i}(\alpha)+\gamma_i)h_i(B(\alpha)) -\sum_{i=1}^k
(L_{i}(\beta-z)+\gamma_i)h_i(B(\beta-z))\\
&=\sum_{i=1}^k L_{i}(\alpha-\beta+z)h_i(y). \ \ \ \ \ \ \ \ \ \ \ \
\ \ \ \ \ \ \ \ \ \ \ \ \ \ \ \ \ \ \ \ \ \ \ \ \ \ \ \  \ \ \ \ \ \
\ \ \ \ \ \ \ \ \ \ \ \ \ \ \ \ \ \ \ \ \ \ \ \ \ \ \ \ \ (2.6)
\end{align*}
Hence by (2.5) and (2.6), we derive that $F(\alpha)=F(\beta-z)$.
Since $F(x)$ is a permutation polynomial of ${\bf F}_{q^m}$, we
obtain that $\alpha-\beta+z=0$. Namely, $x=0$. Thus (2) is true.
The necessity part is proved.

This completes the proof of Theorem 2.1.
\end{proof}

As a special case of Theorem 2.1, we have the following result.\\

\noindent{\bf Corollary 2.1.} {\it Let $L_{1}(x), L_{2}(x) \in {\bf
F}_{q}[x]$ be linearized polynomials. Let $h(x)\in {\bf F}_{q}[x]$
and $\gamma\in {\bf F}_{q^m}$. Then
$F(x):=L_{1}(x)+(L_{2}(x)+\gamma)h(\textup{Tr}_{{\bf F}_{q^m}/{\bf F}_{q}}(x))$
is a
permutation polynomial over ${\bf F}_{q^m}$ if and only if each of the following is true:\\
{\rm (1).} $L_{1}(x)+(L_{2}(x)+\textup{Tr}_{{\bf F}_{q^m}/{\bf
F}_{q}}(\gamma))h(x)\in {\bf F}_{q}[x]$
is a permutation polynomial over ${\bf F}_{q}$.\\
{\rm (2).} For any $y\in {\bf F}_{q}$, $L_{1}(x)+L_{2}(x)h(y)=0$ and
$\textup{Tr}_{{\bf F}_{q^m}/{\bf F}_{q}}(x)=0$ with $x\in {\bf F}_{q^m}$ are both true if
and only if $x=0$.}\\

From Corollary 2.1, we derive the following consequences.\\

\noindent{\bf Corollary 2.2.} \cite{[CHM]} {\it Let
$F(x):=L(x)+xh(\textup{Tr}_{{\bf F}_{q^m}/{\bf F}_{q}}(x))$ with $L(x)\in {\bf
F}_{q}[x]$ being a linearized polynomial and $h(x)\in {\bf
F}_{q}[x]$. Then $F(x)$ is a permutation polynomial over
${\bf F}_{q^m}$ if and only if each of the following is true:\\
{\rm (1).} $L(x)+xh(x)$ is a permutation polynomial over ${\bf F}_{q}$.\\
{\rm (2).} For any $y\in {\bf F}_{q}$, we have that $x\in {\bf F}_{q^m}$
satisfies $L(x)+xh(y)=0$  and $\textup{Tr}_{{\bf F}_{q^m}/{\bf F}_{q}}(x)=0$ if
and only if $x=0$.}

\begin{proof}
This corollary follows from Corollary 2.1 by setting  $L_1(x)=L(x)$,
$L_2(x)=x$ and $\gamma=0.$
\end{proof}

\noindent{\bf Corollary 2.3.} \cite{[M]} {\it Let
$L(x)=a_{0}x+a_{1}x^{q}+\cdots+a_{m-1}x^{q^{m-1}}\in {\bf F}_{q}[x]$
be a linearized polynomial which permutes ${\bf F}_{q^m}$. Let $h(x)\in {\bf
F}_{q}[x]$ and $\gamma\in {\bf F}_{q^m}$. Then the polynomial $F(x):=L(x)+\gamma
h(\textup{Tr}_{{\bf F}_{q^m}/{\bf F}_{q}}(x))$ permutes ${\bf F}_{q^m}$ if and only if the
polynomial $(a_{0}+a_{1}+\cdots+a_{m-1})x+\textup{Tr}_{{\bf F}_{q^m}/{\bf
F}_{q}}(\gamma)h(x)$ permutes ${\bf F}_{q}$.}
\begin{proof}
Since $L(x)$ is a permutation of ${\bf F}_{q^m}$, we have that for any $x\in {\bf F}_{q^m}$,
$L(x)=0$ if and only if $x=0$. So by Corollary 2.1 we know that $F(x)$
is a permutation polynomial over ${\bf F}_{q^m}$ if and only if
$L(x)+\textup{Tr}_{{\bf F}_{q^m}/{\bf F}_{q}}(\gamma)h(x)$ is a permutation
polynomial over ${\bf F}_{q}$.

On the other hand, if $x\in {\bf F}_{q}$, we have
$L(x)=(a_{0}+a_{1}+\cdots+a_{m-1})x$. It then follows that $F(x)$ is
a permutation polynomial over ${\bf F}_{q^m}$ if and only if
$(a_{0}+a_{1}+\cdots+a_{m-1})x+\textup{Tr}_{{\bf F}_{q^m}/{\bf
F}_{q}}(\gamma)h(x)$ is a permutation polynomial over ${\bf F}_{q}$
as desired.
\end{proof}

In what follows we give an example to illustrate Corollary 2.1.\\

\noindent{\bf Example 2.1.} Let ${\bf F}_{q^m}={\bf F}_{8^{m}}$ with $m>1$ being
an odd integer. Let $h(x)=x^{3}-ax$, $L_1(x)=a^{2}x$ and
$L_2(x)=x^{2}$, where $a\in {\bf F}_{8}^{\ast}$. Then
$L_1(x)+L_2(x)h(x)=D_{5}(x,a)$, the Dickson polynomial of degree 5
over ${\bf F}_{8}$. Since ${\rm gcd}(5,q^{2}-1)=1$, by Theorem 7.16
of \cite{[LN]} we know that $D_{5}(x,a)$ is a permutation polynomial
over ${\bf F}_{8}$. That is, $L_1(x)+L_2(x)h(x)=x^{5}-ax^{3}+a^{2}x$
is a permutation polynomial over ${\bf F}_{8}$. Let $y\in {\bf
F}_{8}$ be any element and $x\in {\bf F}_{q^m}$ satisfy that
$\textup{Tr}_{{\bf F}_{q^m}/{\bf F}_{8}}(x)=0$ and $L_1(x)+L_2(x)h(y)=0$. If
$h(y)=0$, then $\textup{Tr}_{{\bf F}_{q^m}/{\bf F}_{8}}(x)=0$ and $L_1(x)=0$.
From $L_1(x)=a^{2}x=0$, we derive that $x=0$. If $h(y)\neq0$, it
then follows from $L_1(x)+L_2(x)h(y)=0$ that $x=0$ or
$x=\frac{a^{2}}{y^{3}-ay}\ne 0$. Assume that
$x=\frac{a^{2}}{y^{3}-ay}$. Since $m$ is odd and
$\frac{a^{2}}{y^{3}-ay}\ne 0$, we have
$$\textup{Tr}_{{\bf F}_{q^m}/{\bf F}_{8}}(x)=\textup{Tr}_{{\bf F}_{q^m}/{\bf F}_{8}}
(\frac{a^{2}}{y^{3}-ay})=\frac{ma^{2}}{y^{3}-ay}\ne 0.$$
Thus we conclude that for any $y\in
K$, $\textup{Tr}_{{\bf F}_{q^m}/{\bf F}_{8}}(x)=0$ and $L_1(x)+L_2(x)h(y)=0$ if
and only if $x=0$. Now by Corollary 2.1, we get that
$$L_1(x)+L_2(x)h(\textup{Tr}_{{\bf F}_{q^m}/{\bf F}_{8}}(x))
=a^{2}x+x^{2}(\textup{Tr}_{{\bf F}_{q^m}/{\bf
F}_{8}}(x)^{3}-a\textup{Tr}_{{\bf F}_{q^m}/{\bf F}_{8}}(x))$$
is a permutation polynomial over ${\bf F}_{q^m}$.\\

\section{\bf Permutation polynomials constructed by the elementary symmetric polynomials}

Throughout this section, let $m$ and $j$ be positive integers such that
$1\leq j\leq m-1$. Let $\sigma_{j}(x_{1}, ..., x_{m})$ denote the
$j$-th elementary symmetric polynomial in $m$ variables $x_{1}, ..., x_{m}$. That is, one has
$$\sigma_{j}(x_{1}, ..., x_{m})=\sum_{1\leq i_1<...<i_j\leq m}x_{i_1}...x_{i_j}.$$
Then we can define the polynomial $\lambda_j(x)$ by
$$\lambda_{j}(x):=\sigma_{j}(x,x^{q},\ldots,x^{q^{m-1}})=\sum_{0\leq i_{1}
 <i_{2}<\ldots<i_{j}\leq m-1}x^{q^{i_{1}}+\ldots+q^{i_{j}}}.$$
Marcos \cite{[M]} used the polynomials $\lambda_{2}(x)$ and
$\textup{Tr}_{{\bf F}_{q^m}/{\bf F}_{q}}(x^{2})$ to give two sufficient
conditions for $xh(\lambda_{2}(x))$ and $xh(\textup{Tr}_{{\bf F}_{q^m}/{\bf
F}_{q}}(x^{2}))$ being permutation polynomials.

In this section, we construct two new classes of permutation
polynomials by using the functions $\lambda_{j}(x)$ and
$\textup{Tr}_{{\bf F}_{q^m}/{\bf F}_{q}}(x^{j})$.
We begin with the following two lemmas which will be needed in the sequel.\\

\noindent{\bf Lemma 3.1.} {\it Let $\alpha \in {\bf F}_{q^m}$ and $a\in {\bf
F}_{q}$. Then $\lambda_{j}(x)\in {\bf F}_{q}[x]$,
$\lambda_{j}(\alpha)\in {\bf F}_{q}$,
$\lambda_{j}(\alpha^{q})=\lambda_{j}(\alpha)$ and $\lambda_{j}(a
\alpha)=a^{j}\lambda_{j}(\alpha)$.}\\

\noindent{\bf Lemma 3.2.} {\it For any integer $j$ satisfying that
$1\leq j\leq m-1$ and ${\rm gcd}(j,q-1)=1$, $\lambda_{j}:
{\bf F}_{q^m}\rightarrow {\bf F}_{q}$ is onto.}
\begin{proof}
First we show that there is an $\alpha \in {\bf F}_{q^m}$ such that $\lambda_{j}(\alpha)\neq 0.$
Since $\lambda_{j}(x)$ has at most
$$\deg(\lambda_{j}(x))=q^{m-j}+...+q^{m-1}\le q+\cdots+q^{m-1}=\frac{q^{m}-1}{q-1}-1<q^{m}=|{\bf F}_{q^m}|$$
roots in ${\bf F}_{q^m}$, there exists an element $\alpha \in {\bf F}_{q^m}$ such that
$\lambda_{j}(\alpha)\neq 0.$ Now pick an $\alpha \in {\bf F}_{q^m}$ such that
$a:=\lambda_{j}(\alpha)\ne 0$. In what follows, we show that for any
$b\in {\bf F}_{q}$, we can find an element $\beta\in {\bf F}_{q^m}$ such that
$\lambda_{j}(\beta)=b$.

Since ${\rm gcd}(j,q-1)=1$, by Theorem 7.8 of \cite{[LN]} we know
that $ax^{j}$ is a permutation polynomial over ${\bf F}_{q}$. It
follows that for any given $b\in {\bf F}_{q}$, there exists an
element $c\in {\bf F}_{q}$ such that $b=ac^{j}$. Since
$\lambda_{j}(\alpha)=a$, letting $\beta:=c\alpha$ gives us that
$$\lambda_{j}(\beta)=\lambda_{j}(c\alpha)=c^{j}\lambda_{j}(\alpha)=ac^{j}=b$$
as desired. Thus Lemma 3.2 is proved.
\end{proof}

Using the polynomials $\lambda_{j}(x)$, we can give the following characterization
on permutation polynomials of the form $x h(\lambda_{j}(x))$, which is the second
main result of this paper. \\

\noindent{\bf Theorem 3.1.} {\it Let $m$ and $j$ be positive
integers such that $1\leq j\leq m-1$ and
$\gcd(j, q-1)=1$. Let $h(x)\in {\bf F}_{q}[x]$. Then
$xh(\lambda_{j}(x))$ is a permutation polynomial over ${\bf F}_{q^m}$ if and
only if $h(0)\neq0$ and $xh(x)^{j}$ permutes ${\bf F}_{q}$. }

\begin{proof} Write $F(x):=xh(\lambda_{j}(x))$.
First we show the sufficiency part. Since $xh(x)^{j}$ permutes ${\bf
F}_{q}$, we obtain that $\delta h(\delta)^{j}\neq0$ for $\delta\in {\bf
F}_{q}^*$. We get that $h(\delta)\neq0$ for $\delta\in {\bf
F}_{q}^{\ast}$. Hence $h(\delta)\neq0$ for all $\delta\in {\bf F}_{q}$.

Now we choose two elements $\alpha,\beta \in {\bf F}_{q^m}$ such that $F(\alpha )=F(\beta )$, namely,
$$\alpha h(\lambda_{j}(\alpha))=\beta h(\lambda_{j}(\beta)).\eqno(3.1)$$
Then $\lambda_j(F(\alpha ))=\lambda_j(F(\beta ))$. Using Lemma 3.1, we infer that
$$\lambda_{j}(\alpha)h(\lambda_{j}(\alpha))^{j}=
\lambda_{j}(\beta)h(\lambda_{j}(\beta))^{j}.\eqno(3.2)$$ Since
$xh(x)^{j}$ permutes ${\bf F}_{q}$, (3.2) tells us that
$\lambda_{j}(\alpha)=\lambda_{j}(\beta)$. It then follows from (3.1)
and the fact that $h(\delta)\neq0$ for all $\delta\in {\bf F}_{q}$
that $\alpha=\beta$. Hence $F(x)$ is a permutation polynomial over
${\bf F}_{q^m}$. The sufficiency part is proved.

Let us now show the necessity part. Assume that $F(x)$ is a
permutation polynomial over ${\bf F}_{q^m}$. First we prove that $h(0)\neq0$. By
Lemma 3.2, we know that the mapping $\lambda_{j}$ is onto if
$\textup{gcd}(j,q-1)=1$. For $1\leq j\leq m-1$, one has that
$$\deg\lambda_{j}(x)=q^{m-j}+...+q^{m-1}\le q+\cdots+q^{m-1}.$$
Thus for any $a\in {\bf F}_{q}^{\ast}$, the equation
$\lambda_{j}(x)=a$ has at most $q+\cdots+q^{m-1}$ roots in ${\bf F}_{q^m}$. Then
the equation $\lambda_{j}(x)=0$ has at least
$q^{m}-(q-1)(q+\cdots+q^{m-1})=q$ roots in ${\bf F}_{q^m}$. Hence
$\lambda_{j}(x)=0$ has a nonzero root in ${\bf F}_{q^m}$. We pick $\alpha\in
{\bf F}_{q^m}^{\ast}$ such that $\lambda_{j}(\alpha)=0$. Then $\alpha h(0)
=\alpha h(\lambda_{j}(\alpha))=F(\alpha )$. Since $F(x)$ is a
permutation polynomial over ${\bf F}_{q^m}$ and $\alpha $ is nonzero, we have
$F(\alpha )\ne 0$. That is, $\alpha h(0)\neq0$. Thus $h(0)\neq0$.

It remains to show that $xh(x)^{j}$ permutes ${\bf F}_{q}$. On the
one hand, by Lemma 3.1 we have
$$\lambda_{j}(F(x))=\lambda_{j}(x)h(\lambda_{j}(x))^{j}.\eqno(3.3)$$
In addition, applying Lemma 3.2, we know that for all integer $j$
with $1\le j\le m-1$ and $\textup{gcd}(j,q-1)=1$, $\lambda_{j}(x)$
is a mapping from ${\bf F}_{q^m}$ onto ${\bf F}_{q}$. This implies that
$$
\{xh(x)^{j}:x\in {\bf
F}_{q}\}=\{\lambda_{j}(x)h(\lambda_{j}(x))^{j}:x\in {\bf F}_{q^m}\}. \eqno(3.4)
$$
Since $F(x)$ permutes ${\bf F}_{q^m}$, it then follows from (3.3) and (3.4) that
\begin{eqnarray*}
& &|\{xh(x)^{j}:x\in {\bf F}_{q}\}|\\
&=&|\{\lambda_{j}(x)h(\lambda_{j}(x))^{j}:x\in {\bf F}_{q^m}\}|\\
&=&|\{\lambda_{j}(F(x)): x\in {\bf F}_{q^m}\}|\\
&=&|\{\lambda_{j}(x): x\in {\bf F}_{q^m}\}|=q.
\end{eqnarray*}
Hence $xh(x)^{j}$ permutes ${\bf F}_{q}$. The necessity part is
proved.

The proof of Theorem 3.1 is complete.
\end{proof}

Now define
$\mu_{j}(x):=\sum_{i=0}^{m-1}x^{jq^{i}}=\textup{Tr}_{{\bf F}_{q^m}/{\bf
F}_{q}}(x^{j})$ for $1\leq j\leq q^{m}-1$. Then $\mu_{j}(x)\in {\bf
F}_{q}[x]$, $\mu_{j}(\alpha)\in {\bf F}_{q}$ and $\mu_{j}(a
\alpha)=a^{j}\mu_{j}(\alpha)$ for all $a\in {\bf F}_{q}$ and $
\alpha\in {\bf F}_{q^m}$. Also $\mu_{j}(x)$ is a mapping from ${\bf F}_{q^m}$ onto ${\bf
F}_{q}$ if $\textup{gcd}(j,q^{m}-1)=1$. Replaced $\lambda_{j}(x)$ by
$\mu_{j}(x)$, we can characterize the permutation polynomials of the
form $x h(\mu_{j}(x))$ as follows. Theorem 3.2 is the third main result
of this paper and its proof is similar as that of Theorem 3.1,
and so we just give a sketch of the proof.\\

\noindent{\bf Theorem 3.2.} {\it Let $m$ and $j$ be positive
integers such that $1\leq j\leq q^m-1$ and $\gcd(j, q^m-1)=1$.
Let $h(x)\in {\bf F}_{q}[x]$. Then $xh(\mu_{j}(x))$ is
a permutation polynomial over ${\bf F}_{q^m}$ if and only if
$h(0)\neq0$ and $xh(x)^{j}$ permutes ${\bf F}_{q}$.}

\begin{proof}
We here merely prove that if $x h(\mu_{j}(x))$ is a permutation
polynomial over ${\bf F}_{q^m}$, then $h(0)\neq0$. The other part of the proof
is similar to that of Theorem 3.1.

Assume that $x h(\mu_{j}(x))$ is a permutation polynomial over ${\bf F}_{q^m}$.
Clearly, there exists a nonzero element $\theta$ such that ${\rm
Tr}_{{\bf F}_{q^m}/{\bf F}_{q}}(\theta)=0$. Since $\gcd(j, q^m-1)=1$, $x^j$
permutes ${\bf F}_{q^m}$. So there is a nonzero element $\omega\in {\bf F}_{q^m}$ such that
$\omega^j=\theta$. Therefore ${\rm Tr}_{{\bf F}_{q^m}/{\bf F}_{q}}(\omega^j)=0$,
i.e., $\mu_{j}(\omega)=0$. Then $\omega h(0)=\omega
h(\mu_{j}(\omega))$. Since $x h(\mu_{j}(x))$ is a permutation
polynomial over ${\bf F}_{q^m}$ and $\omega$ is nonzero, we have $\omega
h(0)\neq0$. Thus $h(0)\neq0$. So Theorem 3.2 is proved.
\end{proof}

Picking $j=2$, then the sufficiency part of Theorems 3.1 and 3.2 becomes
Proposition 12 of \cite{[M]}. Evidently, Theorems 3.1 and 3.2 give an answer to the
open problem raised by Zieve in \cite{[Z1]}.

\section{\bf Permutation polynomials constructed by linear translators}

The main idea of this section is to construct permutation
polynomials over finite fields with linear translators. We first recall
the definition of linear translators as follows: \\

\noindent{\bf Definition 4.1.} \cite{[K]} Let $f: {\bf F}_{q^m}\mapsto {\bf
F}_{q}$, $a\in {\bf F}_{q}$ and $\alpha $ be a nonzero element in
${\bf F}_{q^m}$. If $f(x+u\alpha)-f(x)=ua$ for all $x\in {\bf F}_{q^m}$
and $u\in {\bf F}_{q}$, then we say that $\alpha $ is an {\it $a$-linear
translator} of the function $f$. In particular, $a=f(\alpha)-f(0).$\\

Using linear translators to construct permutation polynomials, we are now in a
position to give the fourth main result of this paper.\\

\noindent{\bf Theorem 4.1.} {\it Let $L_1(x)\in {\bf F}_{q^m}[x]$ be a
linearized permutation polynomial of ${\bf F}_{q^m}$ and
$L_{2}(x)\in {\bf F}_{q^m}[x]$ be a linearized polynomial of ${\bf
F}_{q^m}$. Let $b\in {\bf F}_{q}, \gamma \in {\bf F}_{q^m}, \ h:
{\bf F}_{q}\mapsto {\bf F}_{q},\ f: {\bf F}_{q^m}\mapsto {\bf
F}_{q}$ be surjective and $L_1^{-1}L_{2}(\gamma)$ be a $b$-linear
translator of $f$. Then $L_1(x)+L_{2}(\gamma)h(f(x))$ is a
permutation polynomial of ${\bf F}_{q^m}$ if and only if either
$L_{2}(\gamma)=0$ or $x+bh(x)$ is a permutation polynomial of ${\bf
F}_{q}$.}

\begin{proof}
Write $g(x):=x+bh(x)$ and $G(x):=L_1(x)+L_{2}(\gamma)h(f(x))$.

First we show the sufficiency
part. Since $L_1(x)$ is a permutation polynomial over ${\bf F}_{q^m}$,
so is $G(x)$ if $L_{2}(\gamma)=0$. Assume that $L_{2}(\gamma)\neq0$ and
$g(x)$ is a permutation polynomial of ${\bf F}_{q}$. In the following we show
that $G(x)$ is a permutation polynomial of ${\bf F}_{q^m}$. Take any two
elements $x_1, y_1\in {\bf F}_{q^m}$ such that $G(x_1)=G(y_1)$. That is,
$$L_1(x_1)+L_{2}(\gamma)h(f(x_1))=L_1(y_1)+L_{2}(\gamma)h(f(y_1)),\eqno(4.1)$$
which implies that $L_1(x_1-y_1)=aL_{2}(\gamma),$
where $a:=h(f(y_1))-h(f(x_1))\in {\bf F}_{q}$.
But the assumption that $L_1(x)$ is a permutation polynomial over
${\bf F}_{q^m}$ implies that there exists a unique element $\alpha\in
{\bf F}_{q^m}$ such that $L_1(\alpha)=aL_{2}(\gamma)$. Thus
$\alpha =aL_1^{-1}L_{2}(\gamma)$ and $L_1(\alpha)=L_1(x_1-y_1)$.
It follows immediately that $\alpha=x_1-y_1$, i.e.,
$$x_1=y_1+aL_1^{-1}L_{2}(\gamma)). \eqno(4.2)$$
So (4.1) gives us that
$$L_1(aL_1^{-1}L_{2}(\gamma))+L_{2}(\gamma)
h(f(y_1+aL_1^{-1}L_{2}(\gamma)))=L_{2}(\gamma)h(f(y_1)),\eqno(4.3)$$
which is equivalent to
$$aL_2(\gamma)+L_{2}(\gamma)h(f(y_1+aL_1^{-1}L_{2}(\gamma)))
=L_{2}(\gamma)h(f(y_1)).\eqno(4.4)$$
By the assumption, we have $L_{2}(\gamma)\neq0$. So (4.4) is equivalent to
$$a+h(f(y_1+aL_1^{-1}L_{2}(\gamma)))=h(f(y_1)).   \eqno(4.5)$$
Since $L_1^{-1}L_{2}(\gamma)$ is the $b$-linear translator of $f$, one has
$f(y_1+aL_1^{-1}L_{2}(\gamma))-f(y_1)=ab$. Hence (4.5) is equivalent to
$$a+h(f(y_1)+ab)=h(f(y_1)).\eqno(4.6)$$
Clearly (4.6) is equivalent to
$$(f(y_1)+ab)+bh(f(y_1)+ab)=f(y_1)+bh(f(y_1)).$$
In other words, (4.6) is equivalent to
$$g(f(y_1)+ab)=g(f(y_1)).\eqno(4.7)$$

Claim that $a=0$. In fact, if $b=0$, then by (4.6), one has $a=0$ as claimed. If
$b\ne 0$, then it follows from the assumption that $g(x)$ is a permutation polynomial
of ${\bf F}_{q}$ and (4.7) that $a=0$. The claim is proved. Then by the claim and
(4.2), we derive immediately that $x_1=y_1$. This concludes that $G(x)$ is a
permutation polynomial of ${\bf F}_{q^m}$. The sufficiency part is proved.

Now let's prove the necessity part. Let $G(x)$ be a permutation
polynomial of ${\bf F}_{q^m}$. Suppose that $L_{2}(\gamma)\neq0$. In
what follows we show that $g(x)$ is a permutation polynomial of
${\bf F}_{q}$. If $b=0$, then $g(x)=x$, which is, of course, a
permutation polynomial of ${\bf F}_{q}$. If $b\neq 0$, then we choose any
two elements $u_1\in{\bf F}_{q}$ and $u\in {\bf F}_{q}$ such that
$$g(u_1)=g(u_1+bu).\eqno (4.8)$$
Since $f$ is surjective, there exists an element
$v_1\in {\bf F}_{q^m}$ such that $u_1=f(v_1)$. Then (4.8) is equivalent to
$$g(f(v_1))=g(f(v_1)+bu).\eqno(4.9) $$

Replaced $y_1$ and $a$ by $v_1$ and $u$, respectively, then (4.7) becomes (4.9).
Thus the equivalence of (4.3) and (4.7) applied to (4.9) gives us that
$$L_1(v_1)+L_{2}(\gamma)h(f(v_1))=L_1(v_1+uL_1^{-1}L_{2}(\gamma))+L_{2}
(\gamma)h(f(v_1+uL_1^{-1}L_{2}(\gamma))).$$
Namely, $G(v_1)=G(v_1+uL_1^{-1}L_{2}(\gamma))$. But $G(x)$ is a
permutation polynomial of ${\bf F}_{q^m}$. So $v_1=v_1+uL_1^{-1}L_2(\gamma)$.
Since $L_1(x)$ is a permutation polynomial and $L_{2}(\gamma)\neq0$, we
have $L_1^{-1}L_2(\gamma)\ne 0$. Hence $u=0$. Thus $g(x)$ is a permutation
polynomial of ${\bf F}_{q}$. The necessity part is proved.

This completes the proof of Theorem 4.1.
\end{proof}

Letting $L_2(x)=L_1(x)$, Theorem 4.1 gives the main result of Kyureghyan in \cite{[K]}.\\

\noindent{\bf Corollary 4.1.} \cite{[K]} {\it Let $L(x)\in {\bf
F}_{q^m}[x]$ be a linearized permutation polynomial of ${\bf F}_{q^m}$. Let $b\in
{\bf F}_{q}, \gamma \in {\bf F}_{q^m}, \ h: {\bf F}_{q}\mapsto {\bf F}_{q},\
f: {\bf F}_{q^m}\mapsto {\bf F}_{q}$ be surjective and $\gamma$ be a $b$-linear
translator of $f$. Then $L(x)+L(\gamma)h(f(x))$ is a permutation polynomial
of ${\bf F}_{q^m}$ if and only if $x+bh(x)$ is a permutation
polynomial of ${\bf F}_{q}$.}\\

\begin{center}
{\bf Acknowledgements}
\end{center}

The authors would like to thank the anonymous referee and the editor for their
careful reading of the manuscript and helpful comments and corrections.

\end{document}